\documentclass[a4paper,12pt]{article}
\usepackage[english]{babel}
\usepackage[T2A]{fontenc}
\usepackage[cp1251]{inputenc}
\usepackage{amsthm}
\usepackage[tbtags]{amsmath}
\usepackage{amsfonts,amssymb}
\begin{document}

\newtheorem{theorem}{Theorem}
\newtheorem{lemma}{Lemma.}
\newtheorem{proposition}{Proposition.}
\newtheorem{Cor}{Corollary.}

\begin{center}
{\large\bf Local Centrally Essential Subalgebras\\ of Triangular Algebras}
\end{center}
\begin{center}
O.V. Lyubimtsev, A.A. Tuganbaev
\end{center}

\textbf{Key words:} centrally essential algebra, algebra of upper triangular matrices, nilpotence index.

\textbf{Abstract.}
We study local centrally essential subalgebras in the algebra of all upper triangular matrices over a field of characteristic $\ne 2$. It is proved that the algebras of upper triangular $3\times 3$ or $4\times 4$ matrices have only commutative local centrally essential subalgebras. Every algebra of upper triangular matrices of order exceeding $6$ contains a non-commutative local centrally essential subalgebra.

The work of O.V. Lyubimtsev is done under the state assignment No~0729-2020-0055.  A.A. Tuganbaev is supported by Russian Scientific Foundation, project 16-11-10013P.

\textbf{MSC2010 database 16R99; 20K30}

\section{Introduction}

We only consider associative rings which are not necessarily unital.

\textbf{1.1. Centrally essential rings.}\\ A ring $R$ is said to be \textsf{centrally essential} if either it is commutative or for any non-central element $a\in R$, there exist two non-zero central elements $x$, $y$ with $ax=y$. A ring $R$ with non-zero identity element is centrally essential if and only if the following condition $(*)$ holds:\\ for any non-zero element $a\in R$, there exist two non-zero central elements $x,y\in R$ with $ax=y$. Any non-zero ring with zero multiplication is centrally essential but does not satisfy $(*)$.\footnote{Not necessarily unital rings with $(*)$ are considered in \cite{MT18b}.}

Centrally essential rings with non-zero identity element are studied in papers \cite{MT18}, \cite{MT19}, \cite{MT19b}, \cite{MT19c}, \cite{MT20}, \cite{MT20b}, \cite{MT20c}. Every centrally essential semiprime ring with $1\ne 0$ is commutative; see \cite[Proposition 3.3]{MT18}. In \cite{MT18}, examples of non-commutative group algebras over fields are given. For example, if $Q_8$ is the quaternion group of order 8, then its group algebra over the field of order $2$ is a finite local non-commutative centrally essential ring of order $256$. In addition, in \cite{MT19} it is proved that the external algebra of a three-dimensional linear space over the field of order 3 is a finite non-commutative centrally essential ring, as well. In \cite{MT19c}, there is an example of a centrally essential ring $R$ with $1\ne 0$ such that the factor ring of $R$ with respect to the prime radical is not a PI ring. Abelian groups with centrally essential endomorphism rings are considered in \cite{LT}.

The main result of this paper is Theorem 1.2.

\textbf{1.2. Theorem.} For any field $\mathbb{F}$ of characteristic $\ne 2$ and an arbitrary positive integer $n\ge 7$, there exists a local non-commutative centrally essential subalgebra of the algebra $T_n(\mathbb{F})$ of upper triangular $n\times n$ matrices. 

\textbf{1.3. Remark.} For $n\ge 2$ and a field $F$, the complete matrix algebra $M_n(\mathbb{F})$ over $F$ and the algebra $T_n(\mathbb{F})$ of upper triangular matrices over $F$ are not centrally essential, since all idempotents of any centrally essential ring with $1\ne 0$ are central by \cite[Lemma 2.3]{MT18}. 

An algebra $\mathcal{A}$ is said to be \textsf{centrally essential} if $\mathcal{A}$ is a centrally essential ring. In this paper, we consider local centrally essential subalgebras of the algebra $T_n(\mathbb{F})$ of all upper triangular matrices, where $\mathbb{F}$ is a field of characteristic $\ne 2$. In particular, such subalgebras are of interest, since, for $\mathbb{F} = \mathbb{Q}$, they are quasi-endomorphism algebras of strongly indecomposable torsion-free Abelian groups of finite rank $n$. Quasi-endomorphism algebras of all such groups are local matrix subalgebras in algebra $M_n(\mathbb{Q})$ of all matrices of order $n$ over the field $\mathbb{Q}$; e.g., see \cite[Chapter I, \S 5]{KMT03}.
We remark that the algebra $\mathbb{Q}E$ is the quasi-endomorphism algebra of a strongly indecomposable torsion-free Abelian group of prime rank $p$ if and only if $\mathbb{Q}E$ is isomorphic to a local subalgebra of $T_p(\mathbb{\mathbb{Q}})$. Indeed,
$\mathbb{Q}E/J(\mathbb{Q}E)\cong \mathbb{Q}$ in this case; see \cite [Theorem 4.4.12]{Fat}, where $J(\mathbb{Q}E)$ is the Jacobson radical which is nilpotent, since $\mathbb{Q}E$ is Artinian. It follows from the Weddenburn-Malcev theorem\footnote{See, for example, \cite [Theorem 6.2.1]{Drozd94}.}
that $\mathbb{Q}E\cong \mathbb{Q}E_p\bigoplus J(\mathbb{Q}E)$, where $E_p$ is the identity matrix.
It is known that every nilpotent subalgebra of a full matrix algebra $M_n(\mathbb{F})$ over an arbitrary field $\mathbb{F}$ is transformed by conjugation to a nil-triangular subalgebra; see \cite[Chapter 2, Theorem 6]{ST68}. Since diagonal matrices of a local matrix algebra have equal elements on the main diagonal, they are transformed to itself under conjugation. Consequently, the quasi-endomorphism algebras of such Abelian groups can be realized as matrix subalgebras if and only if these subalgebras are conjugated to some local subalgebra of $T_p(\mathbb{\mathbb{Q}})$. The necessary information on Abelian groups can be found in \cite{Fuc15} and \cite{KMT03}.

Let $\mathbb{F}$ be a field and $\mathcal{A}$ a finite-dimensional algebra over $\mathbb{F}$. An element $a$ of the algebra $\mathcal{A}$ is said to be \textsf{nilpotent} if $a^n = 0$ for some $n \in \mathbb{N}$. The minimal value of $n$ with this property is called the \textsf{nilpotence index} of the element $a$. An algebra $\mathcal{A}$ is called a \textsf{nil-algebra} if every element is nilpotent. For a nil-algebra $\mathcal{A}$, the maximal nilpotence index $\nu(\mathcal{A})$ of its elements is called the \textsf{nil-index}. A positive integer $k$ such that $\mathcal{A}^k = (0)$ and $\mathcal{A}^{k-1}\neq (0)$ is called the \textsf{nilpotence index} of the algebra $\mathcal{A}$. If such an integer $k$ exists, then algebra is said to be of \textsf{nilpotent index} $k$.

For a field $\mathbb{F}$ and an associative algebra $\mathcal{A}$ over $\mathbb{F}$, we denote by $Z(\mathcal{A})$ and $J(\mathcal{A})$ the center and the Jacobson radical of the algebra $\mathcal{A}$, respectively. An algebra $\mathcal{A}$ with $1$ is said to be \textsf{local} if the factor algebra $\mathcal{A}/J(\mathcal{A})$ is a division algebra. Further, $\mathcal{A}$ denotes a local subalgebra in the algebra $T_n(\mathbb{F})$ and $N_n(\mathbb{F})$ denotes the subalgebra of nilpotent matrices in $\mathcal{A}$ (i.e., the algebra of strictly upper triangular matrices). We note that any matrix $A\in \mathcal{A}$ is of the form
$$
A = \begin{pmatrix}
\lambda& a_{12} &\ldots & a_{1n}\\
0& \lambda &\ldots & a_{2n}\\
\vdots& \vdots &\ddots & \vdots\\
0& 0 &\ldots & \lambda
\end{pmatrix}.
$$

We denote by $E_{ij}$ the matrix unit, i.e., the matrix with $1$ on the position $(i, j)$ and zeros on the remaining positions; $E_k$ denotes the identity $k\times k$ matrix. We denote by $<S>$ the linear hull of a subset $S$ of some linear space.

\section{Some General Results}

We recall that a ring $R$ is said to be \textsf{semiprime} if $R$ does not contain two non-zero ideals with zero multiplication.
A ring is said to be \textsf{reduced} if it does not contain zero-square elements. The center of a semiprime ring is a reduced ring.

\textbf{2.1. Proposition.} If $R$ is a centrally essential ring and its center is a semiprime ring, then the ring $R$ is commutative.

\textbf{Proof.} This assertion follows from \cite[Theorem 1.3(a)]{MT18b} in the case, where $R$ satisfies condition 1.1$(*)$. 

We consider the general case. Suppose that the ring $R$ is not commutative, i.e., there exist $x, y\in R$ such that $xy - yx\neq 0$.
Since $R$ is a centrally essential ring and $Z(R)$ is a reduced ring, there exist $c, d\in Z(R)$ such that $d = (xy - yx)c\in Z(R)\backslash{0}$. We note that $xd\neq 0$; otherwise $d^2 = (xy - yx)cd = ((xd)y - y(xd))c = 0$, which is impossible. If $xd\notin Z(R)$, then we can repeat the proof of \cite[Theorem 2.1(c)]{MT18b}. Namely, there exists an element $z\in Z(R)$ such that $xdz\in Z(R)\backslash{0}$. We consider the set $I = \{i\in Z(R)\,|\;ix\in Z(R)\}$. It is clear that $dz\in I$. Now assume that $dI = 0$. Then $d(dz) = 0$, $(dz)^2 = 0$ and $dz = 0$; this is a contradiction. Therefore, $di\neq 0$ for some $i\in I$. However, 
$$
di = (xy - yx)ci = ((ix)y - y(ix))c = 0,
$$
and we obtain a contradiction, as well.

We assume that $xd\in Z(R)$. Then $xdy\neq 0$, since otherwise $d^2 = 0$. In addition, $(xy)d = (yx)d$. Therefore, $(xy - yx)d = 0$. However, then we have $d^2 = (xy - yx)cd = 0$; this is a contradiction. Thus, the ring $R$ is commutative.~\hfill$\square$

The following assertion expands Lemma 2.3 of \cite{MT18} to the case of rings which do not necessarily have $1\ne 0$; in addition, $1$ is replaced by an arbitrary idempotent $e$.

\textbf{2.2. Proposition.} In any centrally essential ring, the following condition holds:
$$
\forall n\in \mathbb{N},\, x_1,\ldots, x_n, y_1,\ldots,y_n, r,\, e = e^2\in R,
$$
$$
\begin{cases}
x_1y_1 +\ldots + x_ny_n = e\\
x_1rey_1 +\ldots + x_nrey_n = 0\\
\end{cases} \Rightarrow re = 0. \eqno(1)
$$
In particular, all idempotents of a centrally essential ring are central.

\textbf{Proof.} We assume that $R$ is a centrally essential ring which satisfies condition $(1)$, but $re\neq 0$.
If $re\in Z(R)$, then
$$
re = re^2 = re(x_1y_1 +\ldots + x_ny_n) = x_1rey_1 +\ldots + x_nrey_n = 0;
$$
this is a contradiction.

Let $re\notin Z(R)$. Then there exist two elements $c,d\in Z(R)$ such that $cre = d\neq 0$. 
We note that $d = cre = (cre)e = de$. Therefore,
$$
d = de = d(x_1y_1 +\ldots + x_ny_n) = c(x_1rey_1 +\ldots + x_nrey_n) = 0.
$$
Now let $e = e^2\notin Z(R)$. We have the relations $e\cdot e = e$ and $e(re)e = 0$, where $r = x - ex$ for any $x\in R$. Consequently, $re = 0$ and $xe = exe$. We note that condition (1) remains true if we replace $re$ by $er$. In this case, $ex = exe$. Therefore, all idempotents of the ring $R$ are central.~\hfill$\square$

\textbf{2.3. Proposition.} Let $\mathcal{A}$ be a local subalgebra of $T_n(\mathbb{F})$ with Jacobson radical $J(\mathcal{A})$. 
The algebra $\mathcal{A}$ is centrally essential if and only if $J(\mathcal{A})$ is a centrally essential algebra.

\textbf{Proof.}  Let us have a matrix $A\in J(\mathcal{A})$ with $A\notin Z(J(\mathcal{A}))$. Since $\mathcal{A}$ is a centrally essential algebra, there exists a matrix $B\in Z(\mathcal{A})$ such that $0\neq AB = C\in Z(\mathcal{A})$. Since $J(\mathcal{A})$ is an ideal, $C\in Z(J(\mathcal{A}))$. If $B\notin J(\mathcal{A})$, then $A = CB^{-1}\in Z(J(\mathcal{A}))$; this contradicts to the choice of the matrix $A$.

Conversely, we have a decomposition $\mathcal{A} = \mathbb{F}E_n\bigoplus J(\mathcal{A})$.
Since $\mathbb{F}E_n\subset Z(\mathcal{A})$, 
$$
Z(J(\mathcal{A}))\subset Z(\mathcal{A}). \eqno(2)
$$ 
If $0\neq A\in \mathcal{A}$ and $A\in Z(\mathcal{A})$, then $0\neq AE_n\in Z(\mathcal{A})$. Let $A\notin Z(\mathcal{A})$ and $A\in J(\mathcal{A})$. Then there exists a matrix $B\in Z(J(\mathcal{A}))$ such that $0\neq AB = C\in Z(J(\mathcal{A}))$.
It follows from relation $(2)$ that $B\in Z(\mathcal{A})$ and $C\in Z(\mathcal{A})$.

Let $A\notin J(\mathcal{A})$. Then $A = A' + A''$, where $0\neq A'\in \mathbb{F}E_n$, $A''\in J(\mathcal{A})$. If $A'' = 0$, then 
$A\in Z(\mathcal{A})$. Otherwise, $0\neq A''B\in Z(J(\mathcal{A}))$ for some matrix $B\in Z(J(\mathcal{A}))$. Then
$$
AB = A'B + A''B = BA' + BA'' = BA.
$$
Since $A'B, A''B\in Z(J(\mathcal{A}))$, we have $AB\in Z(J(\mathcal{A}))\subset Z(\mathcal{A})$. We also note that $AB\neq 0$, since the matrix $A$ is invertible.~\hfill$\square$

By considering Remark 1.3, we obtain the following corollary.

\textbf{2.4. Corollary.} The algebra $N_n(\mathbb{F})$ for $n\ge 3$ is not a centrally essential algebra.

It follows from Proposition 2.3 that the problem of constructing local centrally essential subalgebras of the algebra $T_n(\mathbb{F})$ is equivalent to the problem of constructing centrally essential subalgebras of the algebra $N_n(\mathbb{F})$.

Let $\mathcal{A}$ be a subalgebra of the algebra $N_n(\mathbb{F})$ of nilpotence index $n$. We assume that $\nu(\mathcal{A}) = n$. Then there exists a matrix $A\in \mathcal{A}$ such that $A^{n-1}\neq 0$. We transform $A$ to the Jordan form,
$$
A = E_{12} + E_{23} +\ldots +E_{(n-1)n},
$$
and pass to the corresponding conjugated subalgebra $\mathcal{A}_c$. We denote by $C(A)$ the centralizer of the matrix $A$ in $\mathcal{A}_c$. 
Since the minimal polynomial of the matrix $A$ is equal to its characteristic polynomial, 
$C(A) = \mathbb{F}[A]$, where $\mathbb{F}[A]$ is the ring of all matrices which can be represented in the form $f(A)$, $f(x)\in \mathbb{F}[x]$; 
see \cite[Chapter 1, Theorem 5]{ST68}. For $B\in C(A)$, we have
$$
B = f(A) = \alpha_0E_n + \alpha_1A +\ldots +\alpha_{n-1}A^{n-1}. 
$$
In addition, $\alpha_0 = 0$, since the matrix $B$ is nilpotent. 

\textbf{2.5. Lemma.} If $Z(\mathcal{A}_c) = C(A)$, then the algebra $\mathcal{A}_c$ is commutative.

\textbf{Proof.} Indeed, if $A'\notin C(A)$, then $AA'\neq A'A$. 
However, $A\in C(A) = Z(\mathcal{A}_c)$. This is a contradiction.~\hfill$\square$

\textbf{2.6. Lemma.} Let $\mathcal{A}_c$ be a centrally essential algebra and $Z(\mathcal{A}_c) = <A^{n-1}>$. Then the algebra $\mathcal{A}_c$ is commutative.

\textbf{Proof.} Indeed, if $\mathcal{A}_c$ is not commutative, then for the matrix $A'\notin Z(\mathcal{A}_c)$ we have $BA' = 0$ for
any matrix $B\in Z(\mathcal{A}_c)$.~\hfill$\square$

\section{Nilpotent Centrally Essential Subalgebras of Algebras $N_3(\mathbb{F})$ and $N_4(\mathbb{F})$}

In what follows, we assume that the ground field $\mathbb{F}$ is of characteristic $\ne 2$. 

\textbf{3.1. Proposition.} Any centrally essential subalgebra of the algebra $N_3(\mathbb{F})$ is commutative.

\textbf{Proof.} Every matrix $A\in N_3(\mathbb{F})$ is of the form
$$
A = \begin{pmatrix}
0& a & b\\
0& 0 &c\\
0& 0 & 0
\end{pmatrix}.
$$
Let $\mathcal{A}$ be a non-commutative centrally essential subalgebra of the algebra $N_3(\mathbb{F})$ of nilpotence index 3.
Then $\nu(\mathcal{A}) = 3$. Let a matrix $A\in \mathcal{A}$ be of nilpotence index 3. We transform $A$ to the Jordan form:
$A = E_{12} + E_{23}$. Now, if $B\in C(A)$, then $B = \alpha_1A +\alpha_2A^2$. We note that $Z(\mathcal{A}_c)\subseteq C(A)$; in addition, $\nu(Z(\mathcal{A}_c)) = 3$ by Lemma 2.6. 
However, then $Z(\mathcal{A}_c) = C(A)$ and the algebra $\mathcal{A}_c$ is commutative by Lemma 2.5. This is a contradiction.~\hfill$\square$

In particular, it follows from Proposition 3.1 that all centrally essential endomorphism rings of strongly indecomposable Abelian torsion-free groups of rank 3 are commutative; cf. \cite[Example 3.4]{LT}.

\textbf{3.2. Proposition.} Any centrally essential subalgebra $\mathcal{A}$ of the algebra $N_4(\mathbb{F})$ is commutative.

\textbf{Proof.} If the algebra $\mathcal{A}$ is of nilpotence index $2$, then it is commutative. Let the nilpotence index of $\mathcal{A}$ be equal to $4$. 
There exists $A\in \mathcal{A}$ such that $A^3\neq 0$. Indeed, the algebra $\mathcal{A}$ contains three matrices $S = (s_{ij})$, $T = (t_{ij})$, $P = (p_{ij})$ with $s_{12}\neq 0$, $t_{23}\neq 0$, $p_{34}\neq 0$. Otherwise, the nilpotence index $\mathcal{A}$ is lower than $4$. 
As the required matrix, we can take a matrix $A = (a_{ij})$ such that $a_{i(i+1)}\neq 0$, $i = 1, 2, 3$.
We transform $A$ to the Jordan form,
$$
A = E_{12} + E_{23} + E_{34},
$$
and pass to the corresponding conjugated subalgebra $\mathcal{A}_c$. For a matrix $B\in C(A)$, we have
$$
B = \alpha_1A + \alpha_2A^2 + \alpha_3A^3,
$$
where $\alpha_1, \alpha_2, \alpha_3\in \mathbb{F}$.
It follows from Lemma 2.5 and Lemma 2.6, then $Z(\mathcal{A}_c)\neq C(A)$ and $Z(\mathcal{A}_c)\neq <A^3>$ provided the algebra $\mathcal{A}_c$ is not commutative. Then any matrix $C\in Z(\mathcal{A}_c)$ is of the form
$$
C = \begin{pmatrix}
0 & 0 & c_{13} & c_{14}\\
0 & 0 & 0 & c_{13}\\
0 & 0 & 0 &0\\
0 & 0 & 0 & 0\\
\end{pmatrix}.
$$
Since $\mathcal{A}_c$ is a centrally essential algebra, we have that for the non-zero matrix $D\notin Z(\mathcal{A}_c)$, there exists a matrix $C\in Z(\mathcal{A}_c)$ such that $0\neq DC\in Z(\mathcal{A}_c)$. Since the matrix $D$ is nilpotent, $\textbf{tr} D = 0$. In addition, $\mathcal{A}_c$ is local; therefore all elements on the main diagonal of the matrix $D$ are equal to zero. In this case, it is directly calculated that 
$$
D = \begin{pmatrix}
0 & d_{12} & d_{13} & d_{14}\\
0 & 0 & d_{23} & d_{24}\\
0 & 0 & 0 & d_{12}\\
0 & 0 & 0 & 0\\
\end{pmatrix}.
$$ 
If $d_{12} = 0$ and $D\notin Z(\mathcal{A}_c)$, then $DC = 0$ for any matrix $C\in Z(\mathcal{A}_c)$; this is a contradiction. 
Let $d_{12}\neq 0$ and $DF \neq FD$ for some matrix $F = (f_{ij})\in \mathcal{A}_c$. We find an element $\lambda\in \mathbb{F}$ such that $f_{12} = \lambda d_{12}$. We set $G = \lambda D - F$, $G = (g_{ij})$. Then 
$$
FG = F(\lambda D - F) = \lambda FD - F^2,
$$
$$
GF = (\lambda D - F)F = \lambda DF - F^2.
$$
Therefore, $G\notin Z(\mathcal{A}_c)$ and $g_{12} = 0$. It follows from the obtained contradiction that the algebra $\mathcal{A}_c$ is commutative.

Let the nilpotence index of the algebra $\mathcal{A}$ be equal to 3. Then $\nu(\mathcal{A}) = 3$, i.e., $\mathcal{A}$ contains a matrix $A$ such that $A^2\neq 0$. Indeed, let us assume the contrary,
$A^2 = 0$ for all $A\in \mathcal{A}$. If $A\notin Z(\mathcal{A})$, then $0\neq AB\in Z(\mathcal{A})$ for some matrix $B\in Z(\mathcal{A})$.
Then
$$
(A + B)^2 = A^2 + 2AB + B^2 = 2AB = 0.
$$
Hence $AB = 0$. This is a contradiction.

We transform the matrix $A$ to the Jordan form,
$$
A = E_{12} + E_{23}.
$$
In the corresponding conjugated subalgebra $\mathcal{A}_c$, the centralizer $C(A)$ consists of matrices $B$ of the form
$$
B = \begin{pmatrix}
0 & b_{12} & b_{13} & b_{14}\\
0 & 0 & b_{12} & 0\\
0 & 0 & 0 & 0\\
0 & 0 & b_{43} & 0\\
\end{pmatrix}; \eqno(3)
$$ 
see \cite[Chapter 3, \S 1]{ST68}. In addition, if $C\in Z(C(A))$, then we have
$$
C = \begin{pmatrix}
0 & c_{12} & c_{13} & 0\\
0 & 0 & c_{12} & 0\\
0 & 0 & 0 & 0\\
0 & 0 & 0 & 0\\
\end{pmatrix}. \eqno(4)
$$ 
Let $Z(\mathcal{A}_c)$ be of nilpotence index $3$. Then we can take a matrix from $Z(\mathcal{A}_c)$ as the matrix $A$; see 
\cite[Chapter 1, Proposition 5, Corollary]{ST68}.
In this case, all matrices from $\mathcal{A}_c$ are contained in $C(A)$. Then $\mathcal{A}_c$ consists of the matrices of the form (3)
and the matrices from $Z(\mathcal{A}_c)$ is  form~$(4)$.
If $B = (b_{ij})\notin Z(\mathcal{A}_c)$ and $b_{12} = 0$, then $BC = 0$ for all $C\in Z(\mathcal{A}_c)$. Then $\mathcal{A}_c$ is not a
centrally essential algebra.
Let $b_{12}\neq 0$ and $BD\neq DB$ for some matrix $D = (d_{ij})\in\mathcal{A}_c$. Let $d_{12} = \lambda b_{12}$ and $F = \lambda B - D$, 
$F = (f_{ij})$. Then $f_{12} = 0$ and $F\notin Z(\mathcal{A}_c)$. This is a contradiction.

Let $Z(\mathcal{A}_c)$ be of nilpotence index 2. Then for $C\in Z(\mathcal{A}_c)$, we obtain
$$
C = \begin{pmatrix}
0 & 0 & c_{13} & 0\\
0 & 0 & 0 & 0\\
0 & 0 & 0 & 0\\
0 & 0 & 0 & 0\\
\end{pmatrix}. 
$$
It follows from relation $AC = CA$ for $A\in \mathcal{A}_c$ that
$$
A = \begin{pmatrix}
0 & a_{12} & a_{13} & a_{14}\\
0 & 0 & a_{23} & a_{24}\\
0 & 0 & 0 & 0\\
0 & a_{42} & a_{43} & 0\\
\end{pmatrix}. 
$$
However, then $AC = 0$ for any matrix $C\in Z(\mathcal{A}_c)$.
Consequently, if $\mathcal{A}_c$ is a centrally essential algebra, then $\mathcal{A}_c$ is commutative.~\hfill$\square$

\section{The Proof of Theorem 1.2}

In \cite[Proposition 2.5]{MT19}, it is proved that the external algebra $\Lambda(V)$ over a field $\mathbb{F}$ of characteristic $\ne  2$ is a centrally essential algebra if and only if the dimension of the space $V$ is odd. By considering the regular matrix representation of the algebra $\Lambda(V)$,
we obtain that for an odd positive integer $n > 1$, there exists a non-commutative centrally essential subalgebra of the algebra $N_{2^n}(\mathbb{F})$; also see 
\cite[Example 3.5]{LT}. Therefore, the minimal order of matrices of a non-commutative the external centrally essential algebra is equal to 8.
In the next example, we construct a non-commutative centrally essential algebra of $7\times 7$ matrices.

\textbf{4.1. Example.} We consider a subalgebra $\mathcal{A}$ of $N_{7}(\mathbb{F})$ consisting of matrices $A$ of the form 
$$
A = 
\left(\begin{matrix}
0 & a & b & c & d & e & f\\
0 & 0 & 0 & b & 0 & 0 & d\\
0 & 0 & 0 & 0 & 0 & 0 & e\\
0 & 0 & 0 & 0 & 0 & 0 & 0\\
0 & 0 & 0 & 0 & 0 & 0 & a\\
0 & 0 & 0 & 0 & 0 & 0 & b\\
0 & 0 & 0 & 0 & 0 & 0 & 0\\
\end{matrix}\right).
$$
Let $A'\in \mathcal{A}$, $a' = a + 1$ and the remaining components of the matrix $A'$ coincide with the corresponding components of the matrix $A$. 
Then $AA'\neq A'A$ if $a\neq 0$ and $b\neq 0$. 
Therefore, the algebra $\mathcal{A}$ is not commutative. It is easy to see that $Z(\mathcal{A})$ contain matrices $B$ of the form
$$
B = 
\left(\begin{matrix}
0 & 0 & 0 & c & d & e & f\\
0 & 0 & 0 & 0 & 0 & 0 & d\\
0 & 0 & 0 & 0 & 0 & 0 & e\\
0 & 0 & 0 & 0 & 0 & 0 & 0\\
\ldots & \ldots & \ldots & \ldots & \ldots & \ldots & \ldots\\
0 & 0 & 0 & 0 & 0 & 0 & 0\\
\end{matrix}\right).
$$
If $0\neq A\notin Z(\mathcal{A})$, then $0\neq AB\in Z(\mathcal{A})$ for some matrix $B\in Z(\mathcal{A})$.
Consequently, $\mathcal{A}$ is a centrally essential algebra.

\textbf{4.2. The completion of the proof of Theorem 1.2.}\\
In $N_n(\mathbb{F})$, we consider the subalgebra $\mathcal{A}$ matrices $A$ of the form 
$$
A = 
\left(\begin{matrix}
0 & a_{12} & a_{13} & a_{14} & a_{15} & \ldots & a_{1n-2} & a_{1n-1} & a_{1n}\\
0 & 0 & 0 & a_{13} & 0 & \ldots & 0 & 0 & a_{1\;n-2}\\
0 & 0 & 0 & 0 & 0 & \ldots & 0 & 0 & a_{1\;n-1}\\
0 & 0 & 0 & 0 & 0 & \ldots & 0 & 0 & 0\\
\ldots & \ldots & \ldots & \ldots & \ldots & \ldots & \ldots & \ldots &\ldots\\
0 & 0 & 0 & 0 & 0 & \ldots & 0 & 0 & 0\\
0 & 0 & 0 & 0 & 0 & \ldots & 0 & 0 & a_{12}\\
0 & 0 & 0 & 0 & 0 & \ldots & 0 & 0 & a_{13}\\
0 & 0 & 0 & 0 & 0 & \ldots & 0 & 0 & 0\\
\end{matrix}\right).
$$
We remark that the algebra $\mathcal{A}$ is not commutative; also see Example 4.1. If $B\in Z(\mathcal{A})$, then
$$
B = 
\left(\begin{matrix}
0 & 0 & 0 & b_{14} & b_{15} & \ldots & b_{1n-2} & b_{1n-1} & b_{1n}\\
0 & 0 & 0 & 0 & 0 & \ldots & 0 & 0 & b_{1n-2}\\
0 & 0 & 0 & 0 & 0 & \ldots & 0 & 0 & b_{1n-1}\\
0 & 0 & 0 & 0 & 0 & \ldots & 0 & 0 & 0\\
\ldots & \ldots & \ldots & \ldots \ldots & \ldots & \ldots & \ldots &\ldots & \ldots\\
0 & 0 & 0 & 0 & 0 & \ldots & 0 & 0 & 0\\
\end{matrix}\right).
$$
For $A = (a_{ij})\notin Z(\mathcal{A})$, we have $a_{12}\neq 0$, $a_{13}\neq 0$. 
Let $B = (b_{ij})\in Z(\mathcal{A})$ and $b_{1n-2} = a_{12}$, $b_{1n-1} = a_{13}$. Then $0\neq AB\in Z(\mathcal{A})$.
Indeed, let $AB = C = (c_{ij})$, $BA = D = (d_{ij})$. Then $c_{ij} = d_{ij} = 0$ for all $i\neq 1$, $j\neq n$.
In addition, $c_{1n} = d_{1n} = a_{12}^2 + a_{13}^2$.
Therefore, $\mathcal{A}$ is a centrally essential algebra.~\hfill$\square$

It is known (e.g., see \cite{PV83}) that every rational algebra of dimension $n$ can be realized as the quasi-endomorphism ring of torsion-free Abelian group of rank $n$. By considering \cite[Proposition 3.1]{LT}, we obtain the following corollary.

\textbf{4.3. Corollary.} For every positive integer $n > 6$, there exists an Abelian torsion-free group $A(n)$ of rank $n$ such that its endomorphism ring is a non-commutative centrally essential ring.
 
\section{Remarks and Open Questions}

\textbf{5.1.} It follows from Proposition 2.3 that for a centrally essential local ring $R$, the ring $J(R)$ is centrally essential. The converse is not always true. Indeed, let $R$ be the local ring of upper triangular matrices $2\times 2$ over a division ring $D$ which is not a field. Then $J(R)$ is a commutative ring with zero multiplication. However, the ring $R$ is not centrally essential, since the ring $R/J(R)$ is not commutative; see \cite[Proposition 3.3]{MT19}.

\textbf{5.2. Open question.} Is it true that for any centrally essential ring $R$, the Jacobson radical $J(R)$ is a centrally essential ring?

\textbf{5.3. Open question.} In Theorem 1.2, it is proved that for any positive integer $n > 6$, there exists a local non-commutative centrally essential subalgebra of the algebra $T_n(F)$ of nilpotence index $3$. Is it true that there exists a positive integer $n$ such that $T_n(F)$ contains a local non-commutative centrally essential subalgebra of any nilpotence index $2 < k\le n$?

\textbf{5.4. Open question.} Is it true that there exist local non-commutative centrally essential subalgebras of the algebra $T_n(F)$ for $n = 5$ and $n = 6$?

\textbf{5.5. Open question.} Is it true that for every positive integer $k > 2$, there exists a positive integer $n = n(k)$ such that
the algebra $T_n(F)$ contains a local non-commutative centrally essential subalgebra of nilpotence index $k$?

\end{document}